\documentclass{amsart}
\usepackage{amsmath, amssymb}

\input xy
\xyoption{all}

\newtheorem{theorem}{Theorem}[section]
\newtheorem{lemma}[theorem]{Lemma}
\newtheorem{corollary}[theorem]{Corollary}

\newtheorem{proposition}[theorem]{Proposition}
\newtheorem{claim}[theorem]{Claim}

\theoremstyle{definition}
\newtheorem{example}[theorem]{Example}
\newtheorem{remark}[theorem]{Remark}
\newtheorem{definition}[theorem]{Definition}

\def\cA{\mathcal{A}}
\def\cO{\mathcal{O}}
\def\cI{\mathcal{I}}
\def\cF{\mathcal{F}}
\def\cJ{\mathcal{J}}
\def\cG{\mathcal{G}}
\def\cP{\mathcal{P}}
\def\cM{\mathcal{M}}
\def\cH{\mathcal{H}}
\def\cD{\mathcal{D}}

\def\jet{\mathit{Jet}}
\def\Diff{\mathit{Diff}}
\def\Hom{\mathit{Hom}}

\def\bC{\mathbb{C}}
\def\bN{\mathbb{N}}
\def\bP{\mathbb{P}}

\def\id{\operatorname{id}}
\def\ker{\operatorname{ker}}
\def\hom{\operatorname{Hom}}
\def\diag{\operatorname{diag}}
\def\pr{\operatorname{pr}}
\def\mod{\operatorname{mod}}
\def\d{\operatorname{d}}
\def\diff{\operatorname{Diff}}
\def\grass{\operatorname{Grass}}

\begin{document}

\title{Jet Spaces in Complex Analytic Geometry:\\
An Exposition}
\author{Rahim Moosa}

\address{Massachussetts Institute of Technology \\
Department of Mathematics \\
77 Massachusetts Avenue \\
Cambridge, MA 02139-4307 \\
USA}

\email{moosa@math.mit.edu}

\date{12 March 2004} 

\maketitle

In~\cite{pillay} Pillay described the model-theoretic significance of a result in complex analytic geometry due to Campana~\cite{campana} and Fujiki~\cite{fujiki}.
These notes are an exposition of the ``jet space'' constructions that underly the  Campana/Fujiki theorem.
In particular, we discuss infinitesimal neighbourhoods as well as the sheaves of principal parts, jets, and differential operators.
The material is drawn largely from Grothendieck~\cite{groth} and Kantor~\cite{kantor}.
We also describe how these constructions are used by Campana and Fujiki.

\bigskip

\section{Some sheaf-theoretic preliminaries}
\label{sheaf}

Suppose $f:X\to Y$ is a continuous map of topological spaces, $\cF$ is a sheaf on $X$, and $\cG$ is a sheaf on $Y$.
\begin{itemize}
\item
The {\em direct image of $\cF$}, denoted by $f_*\cF$, is the sheaf on $Y$ which assigns to every open set $V\subseteq Y$ the object $\cF(f^{-1}(V))$.
\item
The {\em inverse image of $\cG$}, denoted by $f^{-1}\cG$, is the sheaf on $X$ associated to the pre-sheaf which assigns to every open set $U\subseteq X$ the object
$$\lim_{\longrightarrow}\{\cG(V):V\supseteq f(U)\}.$$
\end{itemize}
Note that if $X=\{y\}$ is a closed point of $Y$, and $f$ is the inclusion map, then $f^{-1}\cG=\cG_y$, the stalk of $\cG$ at $y$.

The functors $f^{-1}$ and $f_*$ have an adjointness property: there are natural maps $f^{-1}f_*\cF\to \cF$ and $\cG\to f_*f^{-1}\cG$, so that $\hom(f^{-1}\cG,\cF)\approx\hom(\cG,f_*\cF)$.

\smallskip

If $(X,\cO_X)$ and $(Y,\cO_Y)$ are ringed spaces, then a morphism is given by a continuous map on topological spaces, $f:X\to Y$, together with a map of sheaves on $Y$, $f^\#:\cO_Y\to f_*\cO_X$.
By adjointness, a morphism can also be given by a map $f^{-1}\cO_Y\to\cO_X$ of sheaves on $X$.

If $\cG$ is a sheaf of $\cO_Y$-modules, then $f^{-1}\cG$ is a sheaf of $f^{-1}\cO_Y$-modules.
On the other hand, $f^{-1}\cO_Y\to\cO_X$ makes $\cO_X$ into a sheaf of $f^{-1}\cO_Y$-algebras.
Hence we can make $f^{-1}\cG$ into a sheaf of $\cO_X$-modules by tensoring with $\cO_X$ over $f^{-1}\cO_Y$.
This is denoted by
$\displaystyle f^*\cG:=f^{-1}\cG\otimes_{f^{-1}\cO_Y}\cO_X$,
and is sometimes referred to as the {\em $\cO_X$-module inverse image of $\cG$}.

In the particular case when $X=\{y\}$ is a closed (reduced) point of a complex analytic space $Y$, and $f$ is the inclusion map, then $f^*\cG$ is the complex vector space $\cG_y\otimes_{\cO_{Y,y}}\mathbb{C}$.

\begin{remark}
\label{inverse-ideal}
Let $f:X\to Y$ be a morphism of ringed spaces.
Suppose $Y'\subseteq Y$ is a closed subspace whose defining ideal sheaf is $\cJ$.
The inclusion map $\cJ\to\cO_Y$ induces an inclusion $f^{-1}\cJ\to f^{-1}\cO_Y$, and hence a canonical $\cO_X$-linear map
$$ f^*\cJ=f^{-1}\cJ\otimes_{f^{-1}\cO_Y}\cO_X\longrightarrow f^{-1}\cO_Y\otimes_{f^{-1}\cO_Y}\cO_X=\cO_X.$$
Note, however, that this map may not by an inclusion, and so $f^*\cJ$ need not be an ideal sheaf of $\cO_X$.
We therefore define the {\em inverse image ideal sheaf of $\cJ$}, denoted by $f^{-1}\cJ\cdot\cO_X$, to be the image of $f^*\cJ$ in $\cO_X$ under the above map.
This notation is justified by the fact that if $x\in X$ and $f(x)=y\in Y'$, then $(f^{-1}\cJ\cdot\cO_X)_x$ is in fact the ideal sheaf generated by the image of $f^{-1}(\cJ_y)$ under the inclusion $f^{-1}\cO_{Y,y}\to\cO_{X,x}$.
Note also that $\cO_X/(f^{-1}\cJ\cdot\cO_X)\approx f^*(\cO_Y/\cJ)$.
Finally, we say that a closed subspace $X'\subseteq X$ is the {\em inverse image of $Y'$} if on underlying spaces $f^{-1}(Y')=X'$ and the defining ideal sheaf of $X'$ is $f^{-1}\cJ\cdot\cO_X$.
\end{remark}

\begin{lemma}
\label{inverse-ideal=cartesian}
Suppose we have the following commuting diagram of morphisms of ringed spaces
$$\xymatrix{
X'\ar[r]^i\ar[d]_g\ar[dr]^p & X\ar[d]^f \\
Y'\ar[r]^j & Y
}$$
where $i$ and $j$ are closed immersions.
Then $i(X')$ is the inverse image of $j(Y')$ if and only if
$$g^{-1}\cO_{Y'}\otimes_{p^{-1}\cO_Y}i^{-1}\cO_X \ \approx \ \cO_{X'}$$
under the canonical map.
\end{lemma}

\begin{proof}
Let $\cI$ be the defining ideal sheaf of $i(X')$ and $\cJ$ the defining ideal sheaf of $j(Y')$.
Then, as $\cO_X/(f^{-1}\cJ\cdot\cO_X)=f^*(\cO_Y/\cJ)$ we have that $\cI=f^{-1}\cJ\cdot\cO_X$ if and only if
$$\cO_X/\cI \ = \ f^{-1}(\cO_Y/\cJ)\otimes_{f^{-1}\cO_Y}\cO_X.$$
Now, applying $i^{-1}$, this is equivalent to:
$$i^{-1}(\cO_X/\cI) \ = \ p^{-1}(\cO_Y/\cJ)\otimes_{p^{-1}\cO_Y} i^{-1}\cO_X$$
Observing that $i^{-1}(\cO_X/\cI)=\cO_{X'}$ and $j^{-1}(\cO_Y/\cJ)=\cO_{Y'}$, this is in turn equivalent to
$\displaystyle \cO_{X'} = g^{-1}\cO_{Y'}\otimes_{p^{-1}\cO_Y}i^{-1}\cO_X$,
as desired.
\end{proof}

\bigskip

\section{Infinitesimal neighbourhoods}
\label{infinitesimal}

The construction we describe here can be found in Grothendieck, Section~$1$ of ~\cite{groth}.

Let $i:Y\to X$ be a closed immersion of complex analytic spaces.
It follows that the induced map of sheaves on $Y$, $i^{-1}\cO_X\to\cO_Y$, is a surjection.
Let $\cI:=\ker(i^{-1}\cO_X\to\cO_Y)$.
For each $n\geq 0$, the {\em $n$th infinitesimal neighbourhood of $Y$ in $X$} is the complex analytic space
$$Y_i^{(n)}:=( Y , i^{-1}\cO_X/\cI^{n+1} ).$$

\begin{remark}
\label{identify-subspace}
If $Y\subseteq X$ is a closed analytic subspace defined by the ideal sheaf $\cJ$ of $\cO_X$, and $i$ is the inclusion map, then we have a natural identification of $Y_i^{(n)}$ with $(Y,\cO_X/\cJ^{n+1})$.
Indeed, note that we view $\cO_X/\cJ^{n+1}$ as a sheaf on $Y$ by restricting this sheaf on $X$ to its support.
In other words, as a sheaf on $Y$, $\cO_X/\cJ^{n+1}=i^{-1}(\cO_X/\cJ^{n+1})$.
We wish to observe that $\cO_{Y_i^{(n)}}$ is isomorphic to $i^{-1}(\cO_X/\cJ^{n+1})$.
By the exactness of $i^{-1}$, we obtain
$$0\longrightarrow i^{-1}(\cJ^{n+1})\longrightarrow i^{-1}\cO_X\longrightarrow i^{-1}(\cO_X/\cJ^{n+1})\longrightarrow 0$$
It suffices therefore to observe that $\cI^{n+1}=i^{-1}(\cJ^{n+1})$.
As $i^{-1}$ commutes with products of ideal sheaves, we need $\cI=i^{-1}\cJ$.
But this follows from the above exact sequence for $n=0$, together with the fact that $\cO_Y=\cO_X/\cJ=i^{-1}(\cO_X/\cJ)$.
\end{remark}

There is a canonical closed immersion $i^{(n)}:Y_i^{(n)}\to X$.
This map is given by $i$ on the underlying spaces and by the quotient map $i^{-1}\cO_X\to\cO_{Y_i^{(n)}}=i^{-1}\cO_X/\cI^{n+1}$ on the structure sheaves.

There is also the closed immersion $i^{(n,0)}:Y\to Y_i^{(n)}$.
It is the identity on underlying spaces and the quotient map $i^{-1}\cO_X/\cI^{n+1}\to\cO_Y$ on structure sheaves.
More generally, for $m\geq n$, there is $i^{(m,n)}:Y_i^{(n)}\to Y_i^{(m)}$ given by the identity on underlying spaces and the quotient maps $i^{-1}\cO_X/\cI^{m+1}\to i^{-1}\cO_X/\cI^{n+1}$ on sheaves.
Note that the following commute:
$$\xymatrix{
Y_i^{(n)}\ar[rr]^{i^{(m,n)}}\ar[dr]_{i^{(n)}} & & Y_i^{(m)} \ar[dl]^{i^{(m)}} \\
& X &
}$$

The construction of the infinitesimal neighbourhood is functorial.
Given
$$\xymatrix{
Y'\ar[r]^{i'}\ar[d]_f & X'\ar[d]^g \\
Y\ar[r]^i & X
}$$
where $i$ and $i'$ are closed immersions, there is a unique morphism $\alpha:Y_{i'}'^{(n)}\to Y_i^{(n)}$ such that the following commutes:
$$\xymatrix{
Y'\ar[r]^{i^{(0,n)}}\ar[d]_f & Y_{i'}'^{(n)}\ar[d]^{\alpha}\ar[r]^{i'^{(n)}} & X'\ar[d]^g \\
Y\ar[r]^{i^{(0,n)}} & Y_i^{(n)}\ar[r]^{i^{(n)}} & X
}$$
Indeed, considering an isomorphic copy of the situation, we may assume $Y\subseteq X$ with defining ideal sheaf $\cJ$, $Y'\subseteq X'$ with defining ideal sheaf $\cJ'$, both $i$ and $i'$ are inclusions, and $f$ is the restriction of $g$ to $Y'$.
As in Remark~\ref{identify-subspace}, $\cO_{Y_i^{(n)}}=\cO_X/\cJ^{n+1}$ and $\cO_{Y_{i'}'^{(n)}}=\cO_{X'}/\cJ'^{n+1}$.
The map $g^{-1}\cO_X\to\cO_{X'}$, given by $g$ on structure sheaves, induces
$$f^{-1}\cO_Y=g^{-1}(\cO_X/\cJ)=g^{-1}\cO_X/g^{-1}\cJ\longrightarrow\cO_{X'}/\cJ'=\cO_{Y'}.$$
It therefore takes $(g^{-1}\cJ)^{n+1}$ to $\cJ'^{n+1}$ and hence induces
$$f^{-1}\cO_{Y_i^{(n)}}=g^{-1}(\cO_X/\cJ^{n+1})=g^{-1}\cO_X/(g^{-1}\cJ)^{n+1}\longrightarrow\cO_{X'}/\cJ'^{n+1}=\cO_{Y_{i'}'^{(n)}}.$$
This map, together with $f$ on the underlying spaces, gives $\alpha:Y_{i'}'^{(n)}\to Y_i^{(n)}$.
Uniqueness follows immediately from the commuting requirement.

\begin{lemma}
\label{cartesian}
Suppose $Y\subseteq X$ and $Y'\subseteq X'$ are closed subspaces of complex analytic spaces, $i:Y\to X$ and $i':Y'\to X'$ are the inclusions, and $Y'$ is the inverse image of $Y$ under a morphism $g:X'\to X$.
Then $Y_{i'}'^{(n)}$ is the inverse image of $Y_i^{(n)}$.
\end{lemma}

\begin{proof}
Let $\cI$ be the defining ideal of $Y$ in $X$ and $\cJ$ the defining ideal of $Y'$ in $X'$.
That $Y'$ is the inverse image of $Y$ under $g$ means that $g^{-1}\cI\cdot\cO_{X'}=\cJ$ (see Remark~\ref{inverse-ideal}).
But as $(g^{-1}\cI\cdot\cO_{X'})^{n+1}=g^{-1}(\cI^{n+1})\cdot\cO_{X'}$, we see that $\cJ^{n+1}$ is the inverse image ideal sheaf of $\cI^{n+1}$.
Now observe that $\cI^{n+1}$ is the defining ideal of $Y_i^{(n)}$ in $X$ and $\cJ^{n+1}$ is the defining ideal of $Y_{i'}'^{(n+1)}$ in $X'$.
\end{proof}

\bigskip

\section{Relative principal parts}

The following construction can be found in Grothendieck, Section~$2$ of ~\cite{groth}.

Let $p:X\to Y$ be a morphism of complex analytic spaces.
The diagonal map $\diag:X\to X\times_Y X$ is a closed immersion into the fibred product.\footnote{We only consider separated topological spaces.}
Let
$$\Delta^{(n)}_{X/Y}:=X_{\diag}^{(n)}$$
be the $n$th infinitesimal neighbourhood of $X$ in $X\times_Y X$ via the diagonal map.
That is, if $\cI$ denotes the kernel of the induced surjection $\diag^{-1}\cO_{X\times_Y X}\to\cO_X$, then $\Delta_{X/Y}^{(n)}=(X,\diag^{-1}\cO_{X\times_Y X}/\cI^{n+1})$.

Letting $D\subseteq X\times_Y X$ be the diagonal, as a closed analytic subspace of $X\times_Y X$ over $Y$, and $\cJ$ the defining ideal sheaf of $D$, we see by Remark~\ref{identify-subspace} that $\Delta_{X/Y}^{(n)}$ is isomorphic to $(D,\cO_{X\times_Y X}/\cJ^{n+1})$.
This latter is the $n$th infinitesimal neighbourhood of the diagonal in $X\times_Y X$, and is what Fujiki uses in~\cite{fujiki}.

We wish to view $\cO_{\Delta_{X/Y}^{(n)}}$ as a sheaf of $\cO_X$-algebras in a natural way.
This is done by viewing $\Delta_{X/Y}^{(n)}$ as an analytic space over $X$ via the following morphism:
$$\xymatrix{
\Delta_{X/Y}^{(n)}\ar[r]^{\diag^{(n)}} & X\times_Y X\ar[r]^{\quad\quad\pr_1} & X
}$$
Indeed, as $\alpha:=\pr_1\circ\diag^{(n)}$ is the identity on the underlying spaces, $\alpha^{-1}\cO_X=\cO_X$ and $\alpha_*\cO_{\Delta_{X/Y}^{(n)}}=\cO_{\Delta_{X/Y}^{(n)}}$.
We thus obtain a map of sheaves $\cO_X\to\cO_{\Delta_{X/Y}^{(n)}}$, which we denote\footnote{This notation is used by Grothendieck~\cite{groth}, and is short for $(\pr_1\circ\diag^{(n)})^\#$.}
by $\pr_1^*$.
Note that at stalks this map is just
$$\cO_{X,x}\stackrel{\pr_1^\#}{\longrightarrow}\cO_{X\times_YX,(x,x)}\to\cO_{X\times_YX,(x,x)}/\cJ_{(x,x)}^{n+1}=\cO_{\Delta_{X/Y}^{(n)},x}$$
where $\cJ$ is the ideal sheaf of the diagonal in $X\times_YX$.

By the {\em principal parts of order $n$ of $X$ relative to $Y$}, denoted by $\cP_{X/Y}^{(n)}$, we mean the structure sheaf of $\Delta_{X/Y}^{(n)}$ viewed as a sheaf of $\cO_X$-algebras via the map $\pr_1^*:\cO_X\to\cO_{\Delta_{X/Y}^{(n)}}$.
That is, $\cP_{X/Y}^{(n)}=\diag^{-1}\cO_{X\times_Y X}/\cI^{n+1}$ together with its sheaf of $\cO_X$-algebras structure described above.

\begin{remark}
\label{universal-diff}
There is also the morphism $\pr_2\circ\diag^{(n)}:\Delta_{X/Y}^{(n)}\to X$, which (being the identity on underlying spaces) produces the map of sheaves on $X$,
$$\d_{X/Y}^{(n)}:=\pr_2^*:\cO_X\to\cP_{X/Y}^{(n)}$$
which is called the {\em universal differential operator of order $n$ for $X$ relative to $Y$}.
Note that this map is {\em not} $\cO_X$-linear.
For example, suppose $Y$ is a point and $X=\bC$.
By Remark~\ref{identify-subspace},
$\Delta_{\bC}^{(1)}\approx(D,\cO_{\bC^2}/\cJ)$, where $\cO_{\bC^2}$ is the sheaf of germs of holomorphic functions on the complex plane, $D\subseteq\bC^2$ is the diagonal, and $\cJ$ is its defining ideal sheaf.
Taking stalks at the origin, $\cO_{\Delta_{\bC}^{(1)},0}\approx\bC[[x,y]]/(x-y)^2$, where $\bC[[x,y]]$ denotes the ring of convergent power series at the origin in $\bC^2$.
It is then not hard to see that $\cP_{\bC,0}^{(1)}$ is isomorphic to $\bC[[x,y]]/(x-y)^2$ together with the $\bC[[u]]$-algebra structure given by $u\mapsto x\mod(x-y)^2$.
On the other hand $\pr_2^*:\cO_{\bC,0}\to\cP_{\bC,0}^{(1)}$ produces the map $\bC[[u]]\to\bC[[x,y]]/(x-y)^2$ given by $u\mapsto y\mod(x-y)^2$.
This is not $\bC[[u]]$-linear, since for example
$$\pr_2^*(u\cdot u)=y^2\mod(x-y)^2\neq xy\mod(x-y)^2=u\cdot\pr_2^*(u).$$
\end{remark}

\begin{lemma}
\label{sum}
As sheaves of $\cO_X$-modules, $\cP_{X/Y}^{(n)}\approx\cO_X\oplus\cI/\cI^{n+1}$.
\end{lemma}

\begin{proof}
Some words of explanation.
Note that each of $\cP_{X/Y}^{(n)}$, $\cO_X$, and $\cI/\cI^{n+1}$ is equipped with a natural $\cO_X$-module structure.
Indeed, on $\cP_{X/Y}^{(n)}$ it is given by what we have been writing as $\pr_1^*:\cO_X\to\cP_{X/Y}^{(n)}$, and on $\cI/\cI^{n+1}$ the $\cO_X$-module structure is inherited from being an ideal sheaf of $\cP_{X/Y}^{(n)}$.

In order to prove the Lemma, it suffice to show that $\pr_1^*:\cO_X\to\cP_{X/Y}^{(n)}$ is a section to the quotient map $\cP_{X/Y}^{(n)}\longrightarrow\cO_X\approx\diag^{-1}(\cO_{X\times_Y X})/\cI$.
Indeed, this would imply that the following exact sequence of $\cP_{X/Y}^{(n)}$-modules
$$0\longrightarrow\cI/\cI^{n+1}\longrightarrow\cP_{X/Y}^{(n)}\longrightarrow\cO_X\longrightarrow 0$$
is in fact a split exact sequence of $\cO_X$-modules, hence $\cP_{X/Y}^{(n)}\approx\cO_X\oplus\cI/\cI^{n+1}$ as desired.

Recall that $\pr_1^*$ is the map on structure sheaves coming from the morphism $\pr_1\circ\diag^{(n)}:\Delta_{X/Y}^{(n)}\to X$.
On the other hand, the quotient map $\cP_{X/Y}^{(n)}\longrightarrow\cO_X$ is the map on structure sheaves determined by the morphism $\diag^{(n,0)}:X\to\Delta_{X/Y}^{(n)}$.
That $\pr_1^*:\cO_X\to\cP_{X/Y}^{(n)}$ is a section to $\cP_{X/Y}^{(n)}\longrightarrow\cO_X$ then follows from the fact that $\diag$ is a section $\pr_1$ and from following commuting diagram (Section~\ref{infinitesimal}):
$$\xymatrix{
X\ar[rr]^{\diag^{(n,0)}}\ar[dr]_{\diag} & & \Delta_{X/Y}^{(n)}\ar[dl]^{\diag^{(n)}} \\
& X\times_Y X
}$$
\end{proof}

The isomorphism of Lemma~\ref{sum} can also be viewed as one of $\cO_X$-algebras, where $\cO_X\oplus\cI/\cI^{n+1}$ has the natural algebra structure in which $\cI/\cI^{n+1}$ is a nilpotent ideal sheaf of order $n+1$.
This algebra structure is given by: multiplication on $\cO_X$ as a sheaf of rings, multiplication on $\cI/\cI^{n+1}$ as an ideal sheaf of $\cP_{X/Y}^{(n)}$, and multiplication between the two by the $\cO_X$-module structure on $\cI/\cI^{n+1}$.

The sheaf of $\cO_X$-modules $\cI/\cI^{n+1}$ is denoted by $\underline{\Omega}_{X/Y}^n$ and called the {\em sheaf of relative $n$-differentials of $X$ over $Y$}.

Note that $\cP_{X/Y}^{(n)}$ is a coherent sheaf of $\cO_X$-algebras of finite type.
Indeed, this follows from the fact that both $\cO_X$ and $\underline{\Omega}_{X/Y}^n$ are coherent sheaves of $\cO_X$-modules of finite type.

\smallskip

We conclude with the following concrete description of the stalks of $\cP_{X/Y}^{(n)}$.

\begin{proposition}
\label{principal-stalk}
Suppose $x\in X$ and $p(x)=y\in Y$.
There is a canonical isomorphism of $\bC$-vector spaces,
$$\cP_{X/Y,x}^{(n)}\otimes_{\cO_{X,x}}\bC\stackrel{\approx}{\longrightarrow} \cO_{X_y,x}/\cM_{X_y,x}^{n+1}$$
where $X_y$ is the fibre of $p:X\to Y$ above $y$, and $\cM_{X_y,x}$ is its maximal ideal at $x$.
In particular, when $Y$ is a point we have
$$\cP_{X,x}^{(n)}\otimes_{\cO_{X,x}}\bC\approx \cO_{X,x}/\cM_{X,x}^{n+1}$$
\end{proposition}

\begin{proof}
This is proved in ~\cite{groth} in full generality.
We consider only the case when $Y$ is a point.
Let $i:Y\to X$ be the section to $p:X\to Y$, given by $i(Y)=x$.
Consider the following commuting square:
$$\xymatrix{
Y\ar[r]^{i}\ar[d]_{i} & X\ar[d]^{\diag} \\
X\ar[r]_{(ip,\id)\quad} & X\times X
}$$
Under the vertical arrows, we can identify $Y$ and $X$ with closed subspaces of $X$ and $X\times X$ respectively.
As such, $Y$ is the inverse image of $X$ under $(ip,\id)$ (see Remark~\ref{inverse-ideal}).
Taking infinitesimal neighbourhoods with respect to the vertical arrows we obtain (by functoriality)
$$\xymatrix{
Y_i^{(n)}\ar[r]^\alpha\ar[d]_{i^{(n)}} & \Delta_X^{(n)}\ar[d]^{\diag^{(n)}}\\
X\ar[r]_{(ip,\id)\quad} & X\times X 
}$$
Moreover, it follows by Lemma~\ref{cartesian}, that $Y_i^{(n)}$ is the inverse image of $\Delta_X^{(n)}$ under $(ip,\id)$.
It follows that $Y_i^{(n)}$ is also the inverse image of $X$ under $\diag^{(n)}$ -- where we use the horizontal arrows (which are also closed immersions) to identify $Y_i^{(n)}$ and $X$ with closed subspaces of $\Delta_X^{(n)}$ and $X\times X$ respectively.
Indeed, the characterisation given by Lemma~\ref{inverse-ideal=cartesian} shows that being the inverse image is independent of the orientation we consider.

On the other hand, considering
$$\xymatrix{
X\ar[r]^{(ip,\id)\quad}\ar[d]_p &  X\times X\ar[d]^{\pr_1} \\
Y\ar[r]_i & X
}$$
we see that $X$ is itself the inverse image of $Y$ under $\pr_1$ -- using the horizontal arrows for the appropriate identifications.
Putting these together yields
$$\xymatrix{
Y_i^{(n)}\ar[r]^{\alpha}\ar[d]_{p\circ i^{(n)}} & \Delta_X^{(n)}\ar[d]^{\pr_1\circ\diag^{(n)}} \\
Y\ar[r]_i & X
}$$
and $Y_i^{(n)}$ is the inverse image of $Y$ under $\pr_1\circ\diag^{(n)}$ -- where we use the horizontal arrows to identify $Y_i^{(n)}$ and $Y$ with closed subspaces of $\Delta_X^{(n)}$ and $X$ respectively.

By Lemma~\ref{inverse-ideal=cartesian}, and noting that $\pr_1\circ\diag^{(n)}$ is what gives the structure sheaf of $\Delta_X^{(n)}$ the $\cO_X$-module structure of $\cP_X^{(n)}$, we get that
$$(p\circ i^{(n)})^{-1}\cO_Y\otimes_{(i\circ p\circ i^{(n)})^{-1}\cO_X}\alpha^{-1}\cP_X^{(n)} \ \ \approx \ \ \cO_{Y_i^{(n)}}$$
As $p\circ i^{(n)}$ is the identity on underlying spaces and $\alpha$ and $i$ agree on underlying spaces, this yields
$$\cO_Y\otimes_{i^{-1}\cO_X}i^{-1}\cP_X^{(n)} \ \ \approx \ \ \cO_{Y_i^{(n)}}$$
But as $Y$ is a point and $i:Y\to X$ is the section $i(Y)=x$,
we have the following degenerations:
$i^{-1}\cP_X^{(n)}=\cP_{X,x}^{(n)}$,
$i^{-1}\cO_X=\cO_{X,x}$,
$\cO_Y=\bC$,
and
$$\cO_{Y_i^{(n)}}=i^{-1}\cO_X/[\ker(i^{-1}\cO_X\to\cO_Y)]^{n+1}=\cO_{X,x}/\cM_{X,x}^{n+1}$$
Plugging these into the isomorphism displayed above proves the Proposition.
\end{proof}

\begin{remark}
\label{fibre}
Let $\Delta_X^{(n)}\to X$ be the map $\pr_1\circ\diag^{(n)}$.
Then Proposition~\ref{principal-stalk} tells us that for $a\in X$, the (sheaf-theoretic) fibre of $\Delta_X^{(n)}$ over $a$ is (canonically isomorphic to) the analytic space $(\{a\},\cO_{X,a}/\cM_{X,a}^{n+1})$.
\end{remark}

\begin{remark}
\label{section}
Suppose $i:Y\to X$ is a section to $p:X\to Y$ that is also a closed immersion.
The structure sheaf of the $n$th infinitesimal neighbourhood of $Y$ in $X$, $\cO_{Y_i^{(n)}}$, is given the structure of an $\cO_Y$-algebra by $p\circ i^{(n)}:Y_i^{(n)}\to Y$.
The argument for Proposition~\ref{principal-stalk} given above actually proves:
{\em There is a canonical $\cO_Y$-linear isomorphism $i^*\cP_{X/Y}^{(n)}\stackrel{\approx}{\longrightarrow}\cO_{Y_i^{(n)}}$.}
\end{remark}

\begin{corollary}
\label{differential-stalk}
Suppose $x\in X$ and $p(x)=y\in Y$.
There is a canonical isomorphism of $\bC$-vector spaces,
$$\underline{\Omega}_{X/Y,x}^{(n)}\otimes_{\cO_{X,x}}\bC\stackrel{\approx}{\longrightarrow} \cM_{X_y,x}/\cM_{X_y,x}^{n+1}.$$
In particular, when $Y$ is a point we have
$$\underline{\Omega}_{X,x}^{(n)}\otimes_{\cO_{X,x}}\bC\approx \cM_{X,x}/\cM_{X,x}^{n+1}$$
\end{corollary}

\begin{proof}
This is an immediate consequence of Lemma~\ref{sum} and Proposition~\ref{principal-stalk}.
\end{proof}

\bigskip

\section{Jets}

The material in this section is drawn from Kantor~\cite{kantor}.

Let $X$ be a complex analytic space, and $\pi:\cO_X\otimes_{\bC}\cO_X\to\cO_X$ the surjective map of sheaves on $X$ determined by the presheaf map $\cO_X(U)\otimes_{\bC}\cO_X(U)\to\cO_X(U)$ given by $f\otimes g\mapsto fg$.
Let $\cJ:=\ker\pi$, and define the {\em sheaf of jets on $X$} to be
$$\jet_X^{(n)}:=(\cO_X\otimes_{\bC}\cO_X)/\cJ^{n+1}.$$
We view $\jet_X^{(n)}$ as a sheaf of $\cO_X$-algebras by the composition
$$\xymatrix{
\cO_X\ar[r]^{p_1\quad\quad} & \cO_X\otimes_{\bC}\cO_X\ar[r] & (\cO_X\otimes_{\bC}\cO_X)/\cJ^{n+1}=\jet_X^{(n)}
}$$
where the second arrow is the quotient map and the first is determined by the presheaf map $\cO_X(U)\to\cO_X(U)\otimes_{\bC}\cO_X(U)$ given by $f\mapsto f\otimes 1$.

\begin{remark}
\label{universal-diffjet}
We let $d_X^{(n)}:\cO_X\to\jet_X^{(n)}$ be the composition of $p_2:\cO_X\to\cO_X\otimes_{\bC}\cO_X$ with the quotient map, where $p_2$ is determined by $f\mapsto 1\otimes f$.
This should not be confused with $\d_X^{(n)}:\cO_X\to\cP_X^{(n)}$, the ``universal differential operator'' discussed in Remark~\ref{universal-diff}.
\end{remark}

\begin{proposition}
\label{jet-principal}
Suppose $(X,\cO_X)$ is a reduced complex analytic space.
Then there is a canonical $\cO_X$-linear map, $\jet_X^{(n)}\to\cP_X^{(n)}$, such that 
$$\xymatrix{
\jet_X^{(n)}\ar[rr] && \cP_X^{(n)} \\
& \cO_X\ar[ul]^{d_X^{(n)}}\ar[ur]_{\d_X^{(n)}}
}$$
commutes.
Moreover, if $X$ is smooth, then $\jet_X^{(n)}\to\cP_X^{(n)}$ is an isomorphism.
\end{proposition}

\begin{proof}
Fix an open set $U\subseteq X$.
As $X$ is reduced we may interpret the elements of $\cO_X(U)$ as (holomorphic) $\bC$-valued functions on $U$.
Consider the map
$$\alpha_U:\cO_X(U)\otimes_{\bC}\cO_X(U)\to\cO_{X\times X}(U\times U)$$
given by $f\otimes g\mapsto f(x)g(y)$.
On the other hand, we have the natural map
$$\beta_U:\cO_{X\times X}(U\times U)\longrightarrow \lim_{V\supseteq \diag(U)}\cO_{X\times X}(V)$$
The composition, $\beta_U\circ\alpha_U$, determines a map, $\gamma:\cO_X\otimes_{\bC}\cO_X\to\diag^{-1}\cO_{X\times X}$, of sheaves on $X$, such that the following commutes:
$$\xymatrix{
\cO_X\otimes_{\bC}\cO_X\ar[rr]^\gamma\ar[dr]_{\pi} && \diag^{-1}\cO_{X\times X}\ar[dl] \\
& \cO_X &
}$$

We claim that $\gamma$ is injective.
If $U\subseteq X$ is an open set of {\em smooth} points, then $\alpha_U:\cO_X(U)\otimes_{\bC}\cO_X(U)\to\cO_{X\times X}(U\times U)$ is injective.
In general, as $X$ is reduced, there is a dense open set $V\subseteq U$, such that $V$ is a set of smooth points of $X$.
We have the following commuting diagram
$$\xymatrix{
\cO_X(U)\otimes_{\bC}\cO_X(U)\ar[r]^{\alpha_U}\ar[d] & \cO_{X\times X}(U\times U)\ar[d] \\
\cO_X(V)\otimes_{\bC}\cO_X(V)\ar[r]^{\alpha_V} & \cO_{X\times X}(V\times V)
}$$
where the vertical arrows are the restriction maps.
As $V$ is dense in $U$, these vertical arrows are injective, and $\alpha_V$ is injective by the smoothness of $V$.
Hence for each open set $U\subseteq X$, $\alpha_U$ is injective.
Taking limits we have that for all $p\in X$, $\alpha_p:(\cO_X\otimes_{\bC}\cO_X)_p\to\cO_{X\times X,(p,p)}$ is injective.
On the other hand,
$$\cO_{X\times X,(p,p)} \ = \ \lim_{p\in U}\cO_{X\times X}(U\times U) \ = \ \lim_{p\in U} \lim_{V\supseteq\diag(U)}\cO_{X\times X}(V)=(\diag^{-1}\cO_{X\times X})_p$$
and so $\beta_p:\cO_{X\times X,(p,p)}\to (\diag^{-1}\cO_{X\times X})_p$ is an isomorphism.
It follows that $\gamma_p:(\cO_X\otimes_{\bC}\cO_X)_p\to(\diag^{-1}\cO_{X\times X})_p$ is injective, as claimed.

Note that for $p\in X$ a smooth point, $\gamma_p$ is given by the canonical isomorphism $\cO_{X,p}\otimes_{\bC}\cO_{X,p}\stackrel{\approx}{\longrightarrow}\cO_{X\times X,(p,p)}$.
Hence if $X$ is smooth, $\gamma$ is an isomorphism.

In any case, we have the following commuting diagram of exact sequences:
$$\xymatrix{
0\ar[r] & \cI\ar[r] & \diag^{-1}\cO_{X\times X}\ar[rr] && \cO_X\ar[r] & 0\\
0\ar[r] & \cJ\ar[r] & \cO_X\otimes_{\bC}\cO_X\ar[u]_\gamma\ar[urr]_\pi
}$$
It follows that under the injection $\gamma$, $\cJ=\cI\cap(\cO_X\otimes_{\bC}\cO_X)$.
Hence $\cJ^{n+1}\subseteq \cI^{n+1}\cap(\cO_X\otimes_{\bC}\cO_X)$.
It follows that $\gamma$ induces a (not necessarily injective) map
$$\jet_X^{(n)}=(\cO_X\otimes_{\bC}\cO_X)/\cJ^{n+1}\longrightarrow \diag^{-1}\cO_{X\times X}/\cI^{n+1}=\cP_X^{(n)}$$
which is an isomorphism when $X$ is smooth.

That this is $\cO_X$-linear follows from the following commuting diagram 
$$\xymatrix{
\cO_{X,p}\otimes_{\bC}\cO_{X,p}\ar[rr]_{\gamma_p} && \cO_{X\times X,(p,p)} \\
& \cO_{X,p} \ar[ul]^{(p_1)_p} \ar[ur]_{(\pr_1^*)_p}
}$$
for each $p\in X$, which can be easily checked.

Similarly, the above diagram with $p_1$ and $\pr_1^*$ replaced by $p_2$ and $\pr_2^*$ also can be seen to commute.
It follows that
$$\xymatrix{
\jet_X^{(n)}\ar[rr] && \cP_X^{(n)} \\
& \cO_X\ar[ul]^{d_X^{(n)}}\ar[ur]_{\d_X^{(n)}}
}$$
commutes, as desired.
\end{proof}

Even when $X$ is not a manifold, the above map is close to being an isomorphism between the sheaf of jets and the sheaf of principal parts.
Namely, as we will see in the next section (Corollary~\ref{jetdual=principaldual}), it induces an isomorphism on the ``dual'' sheaves. (I use quotations as these sheaves are not necessarily locally free in the non-smooth case.)
For now however, we can show that this induced map is at least an injection.

\begin{corollary}
\label{jet-principal-dual}
Suppose $(X,\cO_X)$ is a reduced complex analytic space.
Composition with the map $\jet_X^{(n)}\to\cP_X^{(n)}$ given by Proposition~\ref{jet-principal} induces an injective map
$$\Hom_{\cO_X}(\cP_X^{(n)},\cO_X)\to\Hom_{\cO_X}(\jet_X^{(n)},\cO_X)$$
\end{corollary}

\begin{proof}
This is a corollary of the proof of Proposition~\ref{jet-principal}, whose notation we continue to use.
As $\cP_X^{(n)}|_U=\cP_U^{(n)}$, $\jet_X^{(n)}|_U=\jet_U^{(n)}$, and $\cO_X|_U=\cO_U$ for all open $U\subseteq X$, it suffices to work with global sections. Namely we show that $$\hom_{\cO_X}(\cP_X^{(n)},\cO_X)\to\hom_{\cO_X}(\jet_X^{(n)},\cO_X)$$ is injective.

First of all note that $\jet_X^{(n)}\to\cP_X^{(n)}$ was induced by $\gamma:\cO_X\otimes_{\bC}\cO_X\to\diag^{-1}\cO_{X\times X}$ under the surjections $\cO_X\otimes_{\bC}\cO_X\to\jet_X^{(n)}$ and $\diag^{-1}\cO_{X\times X}\to\cP_X^{(n)}$.
It suffices therefore to show that composition with $\gamma$ induces an injection 
$$\hom_{\cO_X}(\diag^{-1}\cO_{X\times X},\cO_X)\to\hom_{\cO_X}(\cO_X\otimes_{\bC}\cO_X,\cO_X)$$
Let $f:\diag^{-1}\cO_{X\times X}\to\cO_X$ be an $\cO_X$-linear map (of sheaves on $X$) such that $f\circ\gamma=0$.
We wish to show that $f=0$.
Note that for any $x\in X$,
$$(\diag^{-1}\cO_{X\times X})_x=\lim_{x\in U}\cO_{X\times X}(U\times U).$$
It suffices therefore to show that for all open $U\subseteq X$,
$$f_1(U):\cO_{X\times X}(U\times U)\to\diag^{-1}\cO_{X\times X}(U)\stackrel{f(U)}{\longrightarrow}\cO_X(U)$$
is the zero map.
Recall that for every smooth point $x\in U$, $\gamma_x$ is an ismorphism, and hence $f_x$ is the zero map.
It follows that if $U$ is smooth then $f_1(U)=0$.
In general, let $V\subseteq U$ be the dense open set of smooth points.
We have
$$\xymatrix{
\cO_{X\times X}(U\times U)\ar[d] \ar[r]^{\quad\quad f_1(U)} & \cO_X(U) \ar[d] \\
\cO_{X\times X}(V\times V) \ar[r]^{\quad\quad f_1(V)} & \cO_X(V)
}$$
where the vertical arrows are injections as $V$ is dense in $U$.
Since $f_1(V)$ is the zero map, so is $f_1(U)$, as desired.
\end{proof}

\bigskip

\section{Differential operators}

The constructions described here can be found in Kantor~\cite{kantor}.

\subsection{Differential operators in an abstract setting}
Suppose $f:A\to B$ is a homomorphism of rings (commutative and unitary), and $M, N$ are $B$-modules.
We recursively define a class of $A$-linear maps from $M$ to $N$, called {\em differential operators of order $\leq n$} and denoted by $\diff_{B/A}^{(n)}(M,N)$:
\begin{itemize}
\item
$D\in \diff_{B/A}^{(0)}(M,N)$ if $D$ is $B$-linear.
\item
$D\in \diff_{B/A}^{(n+1)}(M,N)$ if for every $b\in B$,  $[D,b]\in\diff_{B/A}^{(n)}(M,N)$, where $[D,b]:M\to N$ is defined by $[D,b](x):=D(bx)-bD(x)$.
\end{itemize}
We view $\diff_{B/A}^{(n)}(M,N)$ as a $B$-module in the natural way: $(bD)(x):=bD(x)$.

\begin{lemma}
\label{iterate}
An $A$-linear map $D:M\to N$ is a differential operator of order $\leq n$ if and only if for every $b_0,\dots,b_n\in B$, $[[\dots[[D,b_0],b_1],\dots],b_k]$ is the zero map.
\end{lemma}

\begin{proof}
A straightforward induction on $n\geq 0$.
\end{proof}

Let $\pi:B\otimes_AB\to B$ be the map given by $x\otimes y\mapsto xy$, and let $J:=\ker\pi$.
View $(B\otimes_AB)/J^{n+1}$ as a $B$-algebra via the map $B\stackrel{p_1}{\longrightarrow}B\otimes_AB\rightarrow (B\otimes_AB)/J^{n+1}$.
We also set $d:B\to(B\otimes_AB)/J^{n+1}$ to be the map $B\stackrel{p_2}{\longrightarrow}B\otimes_AB\rightarrow (B\otimes_AB)/J^{n+1}$, which is {\em not} $B$-linear.
Then $(B\otimes_AB)/J^{n+1}$ represents the functor $\diff_{B/A}^{(n)}$:

\begin{lemma}
\label{represent-diff}
For every $B$-module $N$, composition with $d$ induces a $B$-linear isomorphism
$\displaystyle \hom_B((B\otimes_AB)/J^{n+1},N)\stackrel{\approx}{\longrightarrow}\diff_{B/A}^{(n)}(B,N)$
\end{lemma}

\begin{proof}
We will use $\overline{z}$ to denote the image of $z\in B\otimes_AB $ under the quotient map $B\otimes_AB\to (B\otimes_AB)/J^{n+1}$.

Let $h\in \hom_B((B\otimes_AB)/J^{n+1},N)$, and set $D:=h\circ d:B\to N$.
Note that for any $b,x\in B$,
\begin{eqnarray*}
[D,b](x) & = & D(bx)-bD(x) \\
& = & h(\overline{1\otimes bx})-bh(\overline{1\otimes x}) \\
& = & h(\overline{1\otimes bx}-\overline{b\otimes x}) \\
& = & h((\overline{1\otimes b}-\overline{b\otimes 1})\cdot(\overline{1\otimes x}))
\end{eqnarray*}
Iterating this, we see that for any $b_0,\dots,b_n\in B$,
\begin{eqnarray*}
[[\cdots[[D,b_0],b_1],\dots],b_n](x) & = & h((\overline{1\otimes b_0}-\overline{b_0\otimes 1})\cdots(\overline{1\otimes b_n}-\overline{b_n\otimes 1})\cdot(\overline{1\otimes x})) \\
& = & h(0\cdot(\overline{1\otimes x}))=0
\end{eqnarray*}
for all $x\in B$.
Indeed, just note that each $(1\otimes b_i-b_i\otimes 1)\in J$.
By Lemma~\ref{iterate}, we have that $D\in\diff_{B/A}^{(n)}(B,N)$.

For injectivity, suppose $h_1,h_2\in\hom_B((B\otimes_AB)/J^{n+1},N)$ and set $D_1:=h_1\circ d$ and $D_2:=h_2\circ d$.
If $D_1=D_2$, then for all $b\in B$, $h_1(\overline{1\otimes b})=h_2(\overline{1\otimes b})$.
Also as $h_1,h_2$ are $B$-linear, we have that
$$h_1(\overline{x\otimes y})=h_1((\overline{x\otimes 1})\cdot(\overline{1\otimes y}))=xh_1(\overline{1\otimes y})=xh_2(\overline{1\otimes y})=h_2(\overline{x\otimes y})$$
for all $x,y\in B$.
Since $(B\otimes_AB)/J^{n+1}$ is generated over $A$ by elements of the form $\overline{x\otimes y}$, it follows that $h_1=h_2$, as desired.

For surjectivity, let $D\in\diff_{B/a}^{(n)}(B,N)$.
Set $D_1:B\otimes_AB\to N$ to be the $B$-linear map given by $x\otimes y\mapsto xD(y)$.
Note that $D_1(xx'\otimes y)=xD_1(x'\otimes y)$ for all $x,x',y\in B$.
Also, for any $x,y,b\in B$,
\begin{eqnarray*}
D_1((1\otimes b-b\otimes 1)(x\otimes y)) & = & D_1(x\otimes by)-D_1(bx\otimes y) \\
& = & xD(by)-bxD(y) \\
& = & x[D,b](y)
\end{eqnarray*}
It follows by iteration that for any $b_0,\dots,b_n\in B$
$$D_1((1\otimes b_0-b_0\otimes 1)\cdots(1\otimes b_n-b_n\otimes 1)(x\otimes y))=x[[\dots[[D,b_0],b_1],\dots],b_n](y)$$
for all $x,y\in B$.
By Lemma~\ref{iterate}, and the fact that $D$ is a differential operator of order $\leq n$, we see that $D_1$ vanishes on the ideal of $B\otimes_{A} B$ generated by all elements of the form $(1\otimes b_0-b_0\otimes 1)\cdots(1\otimes b_n-b_n\otimes 1)$.
This ideal is exactly $J^{n+1}$.
Hence $D_1$ induces a $B$-linear map $h:(B\otimes_{\bC}B)/J^{n+1}\to B$.
On the other hand, as $D_1(1\otimes x)=D(x)$ for all $x\in B$, we have that $D=h\circ d$.

This bijection $\hom_B((B\otimes_AB)/J^{n+1},N)\to\diff_{B/A}^{(n)}(B,N)$ is clearly additive.
Moreover, for all $b,x\in B$, we have
$bh\circ d(x)=bh(\overline{1\otimes x})=b(h\circ d)(x)$.
Hence the bijection is a $B$-linear isomorphism, as desired.
\end{proof}

\begin{example}
\label{algebraic}
Consider $D\in\diff_{\bC[z]/\bC}^{(n)}(\bC[z],\bC[z])$, where $z=(z_1,\dots,z_k)$ is a tuple of indeterminates.
Then $D$ can can be written uniquely as $\displaystyle D:=\sum_{|\alpha|\leq n}a_\alpha D^\alpha$, where $\alpha\in\bN^k$, $a_\alpha\in\bC[z]$ and $\displaystyle D^\alpha:=\frac{1}{\alpha!}\prod_{1\leq i\leq k}(\frac{\partial}{\partial z_i})^{\alpha_i}$.
\end{example}
\begin{proof}
Indeed, $\bC[z]\otimes_{\bC}\bC[z]=\bC[z,z']$, and under this identification, the kernel of $\bC[z]\otimes_{\bC}\bC[z]\to\bC[z]$, $J$, is the one generated by elements of the form $z_i'-z_i$.
Letting $u:=(z'-z)\mod J^{n+1}$, we have that $\bC[z,z']/J^{n+1}$ is freely generated over $\bC[z]$ by elements of the form $u^\alpha$, where $\alpha\in\bN^k$ and $|\alpha|\leq n$.
Taking duals, we have that $\hom_{\bC[z]}(\bC[z,z']/J^{n+1},\bC[z])$ is freely generated over $\bC[z]$ by $\{\tilde{D}^\alpha:|\alpha|\leq n\}$, where $\tilde{D}^\alpha(\sum_\beta P_\beta(z)\cdot u_\beta)=P_\alpha(z)$.
So by Lemma~\ref{represent-diff},
$\diff_{\bC[z]/\bC}^{(n)}(\bC[z],\bC[z])$ is freely generated over $\bC[z]$ by $\{\tilde{D}^\alpha\circ d:|\alpha|\leq n\}$, where $d:\bC[z]\to\bC[z,z']/J^{n+1}$ is given by $d(z)=z'\mod J^{n+1}$.
Fixing $\alpha$ we need to show that $\tilde{D}^\alpha\circ d=D^\alpha$ on $\bC[z]$.
Let $\beta\in\bN^k$.
Then $d(z^\beta)=(z+u)^\beta$.
The $u^\alpha$ coefficient of this expression is $\binom{\beta}{\alpha}z^{\beta-\alpha}$.
Hence $\tilde{D}^\alpha(d(z^\beta))=\binom{\beta}{\alpha}z^{\beta-\alpha}$.
That is, as operators on $\bC[z]$, $\tilde{D}^\alpha\circ d=D^\alpha$.
\end{proof}

\begin{example}
\label{analytic}
Consider $D\in \diff_{\cO/\bC}^{(n)}(\cO,\cO)$, where $\cO$ is the ring of germs of holomorphic function at the origin in $\bC^k$.
Fix a system of co-ordinates, $z=(z_1,\dots,z_k)$.
Then $D$ can be written uniquely as $\displaystyle D:=\sum_{|\alpha|\leq n}a_\alpha D^\alpha$, where $a_\alpha\in\cO$ and the $D^\alpha$ are as in Example~\ref{algebraic}.
\end{example}
\begin{proof}
We first consider $\diff_{\bC[z]/\bC}^{(n)}(\bC[z],\cO)$.
By Lemma~\ref{iterate} we have
\begin{eqnarray*}
\diff_{\bC[z]/\bC}^{(n)}(\bC[z],\cO) & = & \hom_{\bC[z]}((\bC[z]\otimes_{\bC}\bC[z])/J^{n+1},\cO) \\
& = & \hom_{\bC[z]}((\bC[z]\otimes_{\bC}\bC[z])/J^{n+1},\bC[z]) \ \otimes_{\bC[z]} \ \cO \\
& = &\diff_{\bC[z]/\bC}^{(n)}(\bC[z],\bC[z]) \ \otimes_{\bC[z]} \ \cO
\end{eqnarray*}
where the second equality uses that $(\bC[z]\otimes_{\bC}\bC[z])/J^{n+1}$ is a free $\bC[z]$-module of finite rank.
It follows from Example~\ref{algebraic} that $\diff_{\bC[z]/\bC}^{(n)}(\bC[z],\cO)$ is generated freely over $\cO$ by the operators $D^\alpha$.
Hence the restriction map
$$\phi:\diff_{\cO/\bC}^{(n)}(\cO,\cO)\to\diff_{\bC[z]/\bC}^{(n)}(\bC[z],\cO)$$
is surjective:
every element of $\diff_{\bC[z]/\bC}^{(n)}(\bC[z],\cO)$ is of the form $\sum_{|\alpha|\leq n}a_\alpha D^\alpha$, where $a_\alpha\in\cO$, and hence is evidently the restriction of an operator from $\diff_{\cO/\bC}^{(n)}(\cO,\cO)$.
Finally, $\phi$ is injective since the restriction map from $\hom_{\bC}(\cO,\cO)$ to $\hom_{\bC}(\bC[z],\cO)$ is injective by a theorem of Krull.
\end{proof}

We can use the above example to understand, to some extent, differential operators on the germs of holomorphic functions on any complex analytic space:

\begin{lemma}
\label{diff-germ}
Let $\cO$ be the ring of germs of holomorphic functions at the origin in $\bC^k$, and let $I\subseteq\cO$ be an ideal that defines the germ of an analytic set $X_\circ$.
\begin{itemize}
\item[(a)]
Suppose $D\in\diff_{\cO/\bC}^{(n)}(\cO,\cO)$ is such that $D(I)\subseteq I$.
Then the induced map $\overline{D}:\cO_{X_\circ}\to\cO_{X_\circ}$ is a differential operator of order $\leq n$.
\item[(b)]
Conversely, if $E\in\diff_{\cO_{X_\circ}/\bC}^{(n)}(\cO_{X_\circ},\cO_{X_\circ})$, then there is $D\in\diff_{\cO/\bC}^{(n)}(\cO,\cO)$ such that $D(I)\subseteq I$ and $E=\overline{D}$.
\end{itemize}
\end{lemma}

\begin{proof}
Let $r:\cO\to\cO/I=\cO_{X_\circ}$ be the quotient map.

For part~$(a)$, as $D(I)\subseteq I$, it induces a $\bC$-linear map $\overline{D}:\cO_{X_\circ}\to\cO_{X_\circ}$.
That this map is in fact a differential operator of order $\leq n$ can be checked rather easily using the recursive definition of $\diff_{\cO_{X_\circ}/\bC}^{(n)}(\cO_{X_\circ},\cO_{X_\circ})$.

We will prove part~$(b)$ by a series of claims.\footnote{As the argument in Kantor seems circular, we give a different proof.}
Fix a local co-ordinate system $z:=(z_1,\dots,z_k)$ for $\cO$.
Let $E\in\diff_{\cO_{X_\circ}/\bC}^{(n)}(\cO_{X_\circ},\cO_{X_\circ})$ and $D\in\diff_{\cO/\bC}^{(n)}(\cO,\cO)$.

{\em Claim 1: If $r(D(z^\alpha))=E(r(z^\alpha))$ for all $\alpha\in\bN^k$ with $|\alpha|\leq n$, then $r(D(f))=E(r(f))$ for all polynomials $f\in \cO$.}
The proof of this claim goes by induction on $n\geq 0$, the case of $n=0$ being clear.
Fixing $1\leq i\leq k$, note that for all $\beta\in\bN^k$ with $|\beta|\leq n-1$,
\begin{eqnarray*}
r([D,z_i](z^\beta)) & = & r(D(z_iz^\beta))-r(z_iD(z^\beta)) \\
& = & E(r(z_iz^\beta))-r(z_i)E(r(z^\beta)) \\
& = & [E,r(z_i)](r(z^\beta))
\end{eqnarray*}
Hence, by induction, $r([D,z_i](f))=[E,r(z_i)](r(f))$ for all polynomials $f\in\cO$.
Now suppose $r(D(g))=E(r(g))$ for {\em some} polynomial $g\in\cO$.
Then,
\begin{eqnarray*}
r(D(z_ig)) & = & r(z_iD(g)+[D,z_i](g)) \\
& = & r(z_i)E(r(g))+[E,r(z_i)](r(g)) \\
& = & E(r(z_ig))
\end{eqnarray*}
This holds for all $i=1,\dots,k$.
By iterating, we conclude that $r(D(f))=E(r(f))$ for all polynomials $f\in \cO$, proving Claim~$1$.

{\em Claim 2: If $r(D(f))=E(r(f))$ for all polynomials $f\in \cO$, then $r\circ D=E\circ r$.}
Let $g\in \cO$.
For each $m>0$, we may write $g=f_m+h_m$ where $f_m$ is a polynomial and $h_m\in\cM^m$, where $\cM$ is the maximal ideal of $\cO$.
Hence,
$$r(D(g))-E(r(g))=r(D(h_m))-E(r(h_m)).$$
Now, it is not hard to see that as $h_m\in\cM^m$, $D(h_m)\in\cM^{m-n}$ -- where we set $\cM^\ell=0$ for $\ell<0$.
So, $r(D(h_m))\in\cM_{X_\circ}^{m-n}$, where $\cM_{X_\circ}$ is the maximal ideal of $\cO_{X_\circ}$.
Similarly, $E(r(h_m))\in\cM_{X_\circ}^{m-n}$, and so $r(D(g))-E(r(g))\in\cM_{X_\circ}^{m-n}$ for all $m>0$.
As $\displaystyle \bigcap_{m>0}\cM_{X_\circ}^{m-n}=0$, we have shown that $r(D(g))=E(r(g))$ for all $g\in\cO$, proving Claim~$2$.

Combining Claims $1$ and $2$, it suffices, in order to prove part~$(b)$, to find $D\in\diff_{\cO/\bC}(\cO,\cO)$ such that $r(D(z^\alpha))=E(r(z^\alpha))$ for all $\alpha\in\bN^k$ with $|\alpha|\leq n$.
For each such $\alpha$, let $a_\alpha\in\cO$ be such that $E(r(z^\alpha))=r(a_\alpha)$.
Now set $\displaystyle D:=\sum_{|\alpha|\leq k}a_\alpha D^\alpha$, with notation as in Example~\ref{analytic}.
This $D$ works.
\end{proof}

\medskip

\subsection{Sheaves of differential operators.}
The following construction is used by Campana in ~\cite{campana} toward the same ends as Fujiki uses principal parts in~\cite{fujiki}.

\begin{definition}
\label{diff}
Suppose $X$ is a complex analytic space.
The presheaf on $X$ given by $U\mapsto\diff_{\cO_X(U)/\bC}^{(n)}(\cO_X(U),\cO_X(U))$ is in fact a sheaf of $\cO_X$-modules; it is a subsheaf of $\Hom_{\bC}(\cO_X,\cO_X)$.
We denote this sheaf by $\Diff_X^{(n)}$ and call it the {\em sheaf of differential operators on $X$ of order $\leq n$}.
\end{definition}

\begin{proposition}
\label{dual-jet}
Suppose $X$ is a complex analytic space.
Composition with the map $d_X^{(n)}:\cO_X\to\jet_X^{(n)}$ from Remark~\ref{universal-diffjet} induces an $\cO_X$-linear isomorphism of sheaves on $X$,
$$\Hom_{\cO_X}(\jet_X^{(n)},\cO_X)\stackrel{\approx}{\longrightarrow}\Diff_X^{(n)}.$$
\end{proposition}

\begin{proof}
Fix an open set $U\subseteq X$.
Let $\cJ(U):=\ker[\cO_X(U)\otimes_{\bC}\cO_X(U)\to\cO_X(U)]$ where the map is the one given by $f\otimes g\mapsto fg$.
By Lemma~\ref{represent-diff} (with $A:=\bC$ and $B:=N:=\cO_X(U)$) we know that composition with
$$\cO_X(U)\stackrel{p_2(U)}{\longrightarrow}\cO_X(U)\otimes_{\bC}\cO_X(U)\longrightarrow\cO_X(U)\otimes_{\bC}\cO_X(U)/\cJ(U)^{n+1}$$
induces an $\cO_X(U)$-linear isomorphism
$$\hom_{\cO_X(U)}(\cO_X(U)\otimes_{\bC}\cO_X(U)/\cJ(U)^{n+1},\cO_X(U))\stackrel{\approx}{\longrightarrow}\diff_{\cO_X(U)/\bC}^{(n)}(\cO_X(U),\cO_X(U)).$$
This determines a presheaf isomorphism $s:\cH\stackrel{\approx}{\longrightarrow}\Diff_X^{(n)}$, where $\cH$ is the presheaf $U\mapsto\hom_{\cO_X(U)}(\cO_X(U)\otimes_{\bC}\cO_X(U)/\cJ(U)^{n+1},\cO_X(U))$.
In particular, as $\Diff_X^{(n)}$ is a sheaf, so is $\cH$.
Now there is a canonical map $\alpha:\Hom_{\cO_X}(\jet_X^{(n)},\cO_X)\to\cH$ which is an isomorphism at stalks.
As $\cH$ is a sheaf, $\alpha$ is an isomorphism.
Hence $s\circ\alpha:\Hom_{\cO_X}(\jet_X^{(n)},\cO_X)\to\Diff_X^{(n)}$ is an isomorphism.
By construction, it is the map induced by composing with $d_X^{(n)}:\cO_X\to\jet_X^{(n)}$.
\end{proof}

\begin{proposition}
\label{dual-principal}
Suppose $X$ is a reduced complex analytic space.
Composition with the map $\d_X^{(n)}:\cO_X\to\cP_X^{(n)}$ from Remark~\ref{universal-diff} induces an $\cO_X$-linear isomorphism of sheaves on $X$,
$$\Hom_{\cO_X}(\cP_X^{(n)},\cO_X)\stackrel{\approx}{\longrightarrow}\Diff_X^{(n)}.$$
\end{proposition}

\begin{proof}
Let $\theta:\Hom_{\cO_X}(\cP_X^{(n)},\cO_X)\to\Hom_{\bC}(\cO_X,\cO_X)$ be the map induced by composition with $\d_X^{(n)}$.
By Proposition~\ref{jet-principal}, this factors as
$$\xymatrix{
\Hom_{\cO_X}(\cP_X^{(n)},\cO_X) \ar[rr] \ar[dr]_\theta && \Hom_{\cO_X}(\jet_X^{(n)},\cO_X)\ar[dl] \\
& \Hom_{\bC}(\cO_X,\cO_X)
}$$
where $\Hom_{\cO_X}(\jet_X^{(n)},\cO_X)\to\Hom_{\bC}(\cO_X,\cO_X)$ is composition with $d_X^{(n)}$.
Moreover, by Corollary~\ref{jet-principal-dual}, $\Hom_{\cO_X}(\cP_X^{(n)},\cO_X)\to\Hom_{\cO_X}(\jet_X^{(n)},\cO_X)$ is injective.
It follows from Proposition~\ref{dual-jet}, that $\theta$ is an injection of $\Hom_{\cO_X}(\cP_X^{(n)},\cO_X)$ in $\Diff_X^{(n)}$.

Fixing $x\in X$, it suffices to show the surjectivity of $\theta$ onto $\Diff_X^{(n)}$, locally at $x$.
Let $\tilde{X}$ be a sufficiently small open neighbourhood of $x$, such that $\tilde{X}$ is (isomorphic to) an analytic set in some open disc $U\subseteq\bC^k$.
We may assume $x$ corresponds to the origin in $\bC^k$.
Let $E\in\Diff_{X,x}^{(n)}(\tilde{X})$.
Taking $\tilde{X}$ and $U$ sufficiently small, we know by Lemma~\ref{diff-germ}, that $E$ is induced by some $D\in\diff_{\cO_k(U)/\bC}^{(n)}(\cO_k(U),\cO_k(U))$, where $\cO_k:=\cO_{\bC^k}$ is the sheaf of germs of holomorphic functions on $\bC^k$.
By smoothness of $U$, we know that $D$ factors through the principal parts
$$\xymatrix{
\cO_k(U) \ar[rr]^D \ar[dr]_{\d_U^{(n)}} && \cO_k(U) \\
& \cP_{\bC^k}^{(n)}(U) \ar[ur]_{D_1}
}$$
where $D_1\in\hom_{\cO_k(U)}(\cP_{\bC^k}^{(n)}(U),\cO_k(U))$.
(Indeed, $\cP_{\bC^k}^{(n)}(U)=\jet_{\bC^k}^{(n)}(U)$ and apply Proposition~\ref{dual-jet}.)
Let $I$ be the kernel of the natural surjection $\cO_k(U)\to\cO_X(\tilde{X})$ -- that is, $I$ is the defining ideal of the analytic set $\tilde{X}$ in $U$.
Note that the kernel of the induced surjection $\cP_{\bC^k}^{(n)}(U)\to\cP_{X}^{(n)}(\tilde{X})$ is generated by $\d_U^{(n)}(I)$ and $\pr_1^*(I)$.
Since $D_1$ is $\cO_k(U)$-linear, it takes $\pr_1^*(I)$ back to $I$.
On the other hand as $D(I)\subseteq I$, $D_1$ also takes $\d_U^{(n)}(I)$ to $I$.
It follows that $D_1$ induces an $\cO_X(\tilde{X})$-linear map $E_1:\cP_{X}^{(n)}(\tilde{X})\to\cO_X(\tilde{X})$, and we have:
$$\xymatrix{
\cO_k(U)\ar[ddd] \ar[rr]^D \ar[dr]_{\d_U^{(n)}} && \cO_k(U) \ar[ddd] \\
& \cP_{\bC^k}^{(n)}(U) \ar[ur]_{D_1} \ar[d] & \\
& \cP_{X}^{(n)}(\tilde{X})\ar[dr]^{E_1} & \\
\cO_X(\tilde{X}) \ar[ur]^{\d_{\tilde{X}}^{(n)}} \ar[rr]_E && \cO_X(\tilde{X})
}$$
That is, $\theta_{\tilde{X}}$ takes $E_1$ to $E$, and we have shown that $\theta$ is locally surjective (and hence surjective) onto $\Diff_X^{(n)}$.
\end{proof}

\begin{corollary}
\label{jetdual=principaldual}
Suppose $(X,\cO_X)$ is a reduced complex analytic space.
Composition with the map $\jet_X^{(n)}\to\cP_X^{(n)}$ given by Proposition~\ref{jet-principal} induces an isomorphism
$$\Hom_{\cO_X}(\cP_X^{(n)},\cO_X)\stackrel{\approx}{\longrightarrow}\Hom_{\cO_X}(\jet_X^{(n)},\cO_X)$$
\end{corollary}

\begin{proof}
This follows from Propositions~\ref{dual-jet} and~\ref{dual-principal}.
\end{proof}

\begin{remark}
Corollary~\ref{jetdual=principaldual} does {\em not} imply that $\jet_X^{(n)}\approx\cP_X^{(n)}$.
If these sheaves were locally free then it would -- by taking duals.
In the case that $X$ is smooth the sheaves of principal parts and of jets are locally free, but in that case we already know $\jet_X^{(n)}\approx\cP_X^{(n)}$ by Proposition~\ref{jet-principal}.
\end{remark}

\bigskip

\section{Applications: the Dichotomoy Theorem}

In this section we describe how the sheaves of principal parts and differential operators are used by Fujiki~\cite{fujiki} and Campana\cite{campana}, respectively, to prove a theorem which Pillay~\cite{pillay} has observed implies the dichotomy theorem for minimal types in the theory of compact complex analytic spaces.

\begin{definition}
\label{proj}
A morphism $f:X\to S$ of (possible non-reduced) irreducible complex analytic spaces is {\em projective} if it factors through a closed embedding into a projective linear space over $S$.
That is, there exists a commutative diagram:
$$\xymatrix{
X\ar[rr]^h\ar[dr]_f && \bP(\cF)\ar[dl] \\
& S}$$
where $h$ is a closed embedding and $\cF$ is a coherent analytic sheaf on $S$.
We say that $g:Y\to S$ is {\em Moishezon} if there is a commutative diagram:
$$\xymatrix{
Y\ar[rr]^\alpha\ar[dr]_g && X\ar[dl]^f \\
& S}$$
where $\alpha$ is a bimeromorphism and $f:X\to Y$ is projective.
\end{definition}

A discussion of the model-theoretic significance of Moishezon morphisms, namely its relationship to {\em internality in $\bP$}, can be found in~\cite{ret}.

\medskip

We work in the following context:
Suppose $X$ and $S$ are reduced and irreducible compact complex analytic spaces, and $Z\subseteq S\times X$ is an irreducible analytic subset.
Let $\rho:Z\to S$ be the first co-ordinate projection map.
We view:
$$\xymatrix{
Z\subseteq S\times X \ar@<-5ex>[d]^\rho \\
S\quad\quad\quad\quad}$$
as a family of analytic subspaces of $X$ parametrised by $S$: namely the family of fibres $Z_s:=\{x\in X:(s,x)\in Z\}$ for $s\in S$.
{\em We make the following standing assumptions}:
\begin{itemize}
\item[(1)]
The map $\rho$ is surjective and there is a dense Zariski open set $W\subseteq S$ such that $Z_w$ is reduced and irreducible for all $w\in W$.
\item[(2)]
$\rho:Z\to S$ is a {\em canonical family}: There is a dense Zariski open set $V\subseteq S$ such that if $v,v'\in V$ and $Z_v=Z_{v'}$, then $v=v'$.
\end{itemize}
In the language of types, these assumptions correspond to saying that we are studying a stationary type over its canonical base.

\smallskip

We aim to expose the following result of Campana~\cite{campana}/Fujiki~\cite{fujiki}:

\begin{theorem}
\label{cf}
The second co-ordinate projection $\pi:Z\to X$ is Moishezon.
\end{theorem}

\medskip

\subsection{Fujiki's approach}

Fix $n\geq 0$ and $a\in X$.
Let $\Delta_X^{(n)}\to X$ be the map $\pr_1\circ\diag^{(n)}$.
Proposition~\ref{principal-stalk} tells us that the (sheaf-theoretic) fibre of $\Delta_X^{(n)}\to X$ over $a$, which we will denote by $\Delta_{X,a}^{(n)}$, is the non-reduced space $(\{a\},\cO_{X,a}/\cM_{x,a}^{n+1})$.
That is, $\Delta_{X,a}^{(n)}$ is the $n$th infinitesimal neighbourhood of $a$ in $X$.

We wish to view $\pi^{-1}(a)\subseteq Z$ as parametrising a family of analytic subsets of $\Delta_{X,a}^{(n)}$.
The basic idea is as follows: for $c\in\pi^{-1}(a)$, $Z_{\rho(c)}$ is an analytic subset of $X$ passing through $a$, and hence the $n$th infinitesimal neighbourhood of $a$ in $Z_{\rho(c)}$ is a subspace of $\Delta_{X,a}^{(n)}$.
We can view $\pi^{-1}(a)$ as parametrising the family of $n$th infinitesimal neighbourhood of $a$ in $Z_{\rho(c)}$, for $c\in\pi^{-1}(a)$.
We will do this formally, and uniformly in $a$, as follows.

Identify $Z\times_SZ$ with a subspace of $Z\times_X(X\times X)$ in the natural way:
$$Z\times_SZ\subseteq Z\times_S(S\times X)\approx Z\times X\approx Z\times_X(X\times X)$$
Also identify $\Delta_X^{(n)}$ with a subspace of $X\times X$ via $\diag^{(n)}$.
Now define
$$Y_n:=(Z\times_SZ)\cap(Z\times_X\Delta_X^{(n)})$$
where the intersection is taking place inside $Z\times_X(X\times X)$.
We have the co-ordinate projection $Z\times_X(X\times X)\to Z$ which restricts to:
$$\xymatrix{
Y_n\subseteq Z\times_X\Delta_X^{(n)}\ar@<-7ex>[d]^{\gamma_n} \\
Z\quad\quad\quad\quad\quad\quad}$$

This makes $Y_n\to Z$ into a family of analytic subspaces of $\Delta_{X}^{(n)}$ over $X$.
Explicitly, for each $a\in X$, restricting everything to $a$ we obtain:
$$\xymatrix{
(Y_n)_{\pi^{-1}(a)}\subseteq \pi^{-1}(a)\times\Delta_{X,a}^{(n)}\ar@<-9ex>[d]^{\gamma_n} \\
\pi^{-1}(a)\quad\quad\quad\quad\quad\quad\quad\quad}$$
and for each $c\in\pi^{-1}(a)\cap W$ the (sheaf-theoretic) fibre $Y_{n,c}$ is the $n$th infinitesimal neighbourhood of $a$ in $Z_{\rho(c)}$ viewed as a subspace of $\Delta_{X,a}^{(n)}$.

For some dense Zariski open set $U_n\subseteq Z$, $\gamma_n:Y_n\to Z$ is flat over $U_n$.
So $U_n$ is a parameter space for a flat family of analytic subspaces of $\Delta_X^{(n)}$ over $X$.
But there is a {\em universal parameter space for flat families of analytic subspaces of $\Delta_X^{(n)}$ over $X$} -- namely the relative Douady space of $\Delta_X^{(n)}$ over $X$.
What this means is that there is a reduced complex analytic space over $X$, $\cD(\Delta_X^{(n)}/X)\to X$, and a diagram of holomorphic maps:
$$\xymatrix{
U_n\ar[rr]^{\tau_n}\ar[dr]_\pi && \cD(\Delta_X^{(n)}/X)\ar[dl]\\
& X}$$
such that for $c,c'\in U_n$, $\tau_n(c)=\tau_n(c')$ if and only if $Y_{n,c}=Y_{n,c'}$.
Moreover, by Hironaka's Flattening Theorem, $\tau_n$ extends to a meromorphic map $Z\to \cD(\Delta_X^{(n)}/X)$, which we will also (abusively) denote as $\tau_n$.

A key fact is that every component $C$ of $\cD(\Delta_X^{(n)}/X)$ is compact and the restricted map $C\to X$ is projective.
This is because $\Delta_X^{(n)}\to X$ is projective (indeed, it is finite), and is proved in Lemma~5 of~\cite{fujiki}.
Hence to show that $\pi$ is Moishezon (i.e. to prove Theorem~\ref{cf}) it will suffice to show that for some $n>0$, $\tau_n$ is a bimeromorphism with its image.

\begin{remark}
Note that $\gamma_n:Y_n\to Z$ does not live definably in the many-sorted structure $\cA$ of compact complex analytic spaces.
This is because $Y_n$ is not reduced, and if we take its reduction then we lose all the information we require -- namely, the infinitesimal neighbourhoods.
However, the map $\tau_n:Z\to \cD(\Delta_X^{(n)}/X)$ is definable and it captures the information given by $\gamma_n:Y_n\to Z$ -- since for $c,c'\in U_n$, $\tau_n(c)=\tau_n(c')$ if and only if $Y_{n,c}=Y_{n,c'}$.
This is somehow the point.
\end{remark}

\smallskip

We now show that for some $n>0$, $\tau_n$ is a bimeromorphism with its image.

First note what it means to prove that some $\tau_n$ is a bimeromorphism with its image.
For $c,c'\in U_n$, supposing $\tau_n(c)=\tau_n(c')$ one asks whether $c=c'$.
As $\pi(c)=\pi(c')$ (since $\tau_n$ is over $X$), this is the same as asking if $\rho(c)=\rho(c')$.
But $\rho:Z\to S$ is a canonical family (standing assumption $\# 2$), so after replacing $U_n$ with $U_n\cap V$, $\rho(c)=\rho(c')$ if and only if $Z_{\rho(c)}=Z_{\rho(c')}$.
On the other hand, from the definition of $\tau_n$ we have that $\tau_n(c)=\tau_n(c')$ if and only if $Y_{n,c}=Y_{n,c'}$; and the latter means that the $n$th infinitesimal neighbourhood of $a:=\pi(c)=\pi(c')$ in $Z_{\rho(c)}$ and in $Z_{\rho(c')}$ agree.
So to say that $\tau_n$ is a bimeromorphism with its image is to say (essentially) that the analytic subsets of $X$ from the family $Z\to S$ are distinguishable by the $n$th infinitesimal neighbourhoods of a common point.

\begin{claim}
\label{claim1}
For irreducible analytic sets $A,B\subseteq X$ with $A\cap B\neq\emptyset$, $A=B$ if and only if at some common point the $n$th infinitesimal neighbourhoods agree for all $n$.
\end{claim}

\begin{proof}
Only the right to left direction requires proof.
We will show that $A\subseteq B$ and concluded by symmetry that $A=B$.
Let $a\in A\cap B$.
As $A\cap B$ is an analytic subset of $A$, and as $A$ is irreducible, it suffices to show that in some non-empty subset of $X$ containing $a$, say $U$, $(A\cap U)\subseteq(B\cap U)$.
This latter condition is local and so we may assume that $A$ and $B$ are analytic subsets of the unit disc in some $\bC^k$, and that $a$ is the origin.
Suppose $f\in \cO$ is a germ of a holomorphic function at the origin which vanishes on the germ of $B$ at the origin.
Hence the image of $f$ in $\cO_B$ -- the ring of germs of holomorphic functions on $B$ at the origin -- is contained in $\cM_B^{n+1}$ for all $n+1$.
As the $n$th infinitesimal neighbourhoods of the origin in A and B agree for all $n$, this implies that the image of $f$ in $\cO_A$ is contained in $\cM_A^{n+1}$ for all $n+1$.
As $\bigcap_n\cM_A^{n+1}=0$, we have that the image of $f$ in $\cO_A$ is zero.
It follows that the defining ideal of the germ of $B$ is contained in that of the germ of $A$.
Hence, for some open set $U$ about the origin, $(A\cap U)\subseteq(B\cap U)$, as desired.
\end{proof}

\begin{claim}
\label{claim2}
Let $m\leq n$ and $c,c'\in U_n\cap U_m$.
If $\tau_n(c)=\tau_n(c')$ then $\tau_m(c)=\tau_m(c')$.
\end{claim}

\begin{proof}
Suppose $\tau_n(c)=\tau_n(c')$ and let $a:=\pi(c)=\pi(c')$.
So the $n$th infinitesimal neighbourhoods of $a$ in $Z_{\rho(c)}$ and $Z_{\rho(c')}$ agree.
It is clear from the definitions that so then do the $m$th infinitesimal neighbourhoods.
Hence $\tau_m(c)=\tau_m(c')$.
\end{proof}

Theorem~\ref{cf} now follows by an elementary model-theoretic argument.

Let $U:=\bigcap_nU_n$ and let $V,W\subseteq S$ be the dense Zariski open sets given by the standing assumptions $1$ and $2$.
Let $\Phi(x,y)$ be the partial type which says that
$$x,y\in U\cap \rho^{-1}(V)\cap\rho^{-1}(W)\text{; }\tau_n(x)=\tau_n(y)\text{ for all }n>0\text{; and }x\neq y.$$
If $(c,c')\models\Phi(x,y)$ then the $n$th infinitesimal neighbourhoods of $a:=\pi(c)=\pi(c')$ in $Z_{\rho(c)}$ and $Z_{\rho(c')}$ agree for all $n$.
These are irreducible analytic subsets of $X$ (as $\rho(c),\rho(c')\in W$), and so $Z_{\rho(c)}=Z_{\rho(c')}$ (by Claim~\ref{claim1}), and so $\rho(c)=\rho(c')$ (as $\rho(c),\rho(c')\in V$), and so $c=c'$.
But this is a contradiction.
Hence $\Phi(x,y)$ is not realised in $Z$.
By $\omega_1$-compactness (note that $\Phi(x,y)$ involves only countably many formulae) there is some $n>0$ such that for all $c,c'\in U_1\cap\cdots\cap U_n\cap \rho^{-1}(V)\cap\rho^{-1}(W)$, if $c\neq c'$ then $\displaystyle \bigvee_{i=1}^n\tau_i(c)\neq\tau_i(c')$.
By Claim~\ref{claim2}, $\tau_n(c)\neq\tau_n(c')$.
So $\tau_n$ is a bimeromorphism with its image, as desired.
\qed

\medskip

\subsection{Campana's approach}

We conclude with a brief sketch of the idea, as we understand it, behind Campana's proof of Theorem~\ref{cf}.

The techniques in~\cite{campana} differ from Fujiki's in essentially two ways.
First of all, the analytic {\em space} of infinitesimal neighbourhoods $\Delta_X^{(n)}$ over $X$ is replaced by the {\em sheaf} of differential operators $\Diff_X^{(n)}$ on $X$.
By Proposition~\ref{dual-principal} this amounts to taking the dual to the structure sheaf of $\Delta_X^{(n)}$.
Secondly, Douady spaces are replaced with Grassmanians.

The Grassmanian of a coherent analytic sheaf on $X$ is a complex analytic space over $X$, introduced by Grothendieck in~\cite{groth2}, which relativises the usual notion of Grassmanian for complex vector spaces.
Fixing $n>0$,
the role of $\cD(\Delta_X^{(n)}/X)\to X$ will now be played by
$\grass_{r_n}(\Diff_X^{(n)})\to X$
for some positive integer $r_n$ (to be determined below).
Fixing $a\in X$, the fibre of $\grass_{r_n}(\Diff_X^{(n)})$ over $a$ can be identified with the Grassmanian of $r_n$-planes (that is, co-dimension $r_n$ subspaces) of the complex vector space $\displaystyle \Diff_{X,a}^{(n)}\otimes_{\cO_{X,a}}\bC$.
For $a\in X$ smooth, the latter can be computed as follows:
\begin{eqnarray*}
\Diff_{X,a}^{(n)}\otimes_{\cO_{X,a}}\bC & = & \Hom_{\cO_X}(\cP_X^{(n)},\cO_X)_a\otimes_{\cO_{X,a}}\bC\\
& = & \hom_{\cO_{X,a}}(\cP_{X,a}^{(n)},\cO_{X,a})\otimes_{\cO_{X,a}}\bC\\
& = & \hom_{\bC}(\cP_{X,a}^{(n)}\otimes_{\cO_{X,a}}\bC \ , \ \bC)\\
& = & \hom_{\bC}(\cO_{X,a}/\cM_{X,a}^{n+1} \ ,\ \bC)
\end{eqnarray*}
where the first equality is by Proposition~\ref{dual-principal}, the second by the choice of $a\in X$ smooth, and the last by Proposition~\ref{principal-stalk}.
Hence
$$\grass_{r_n}(\Diff_X^{(n)})_a=\grass_{r_n}([\cO_{X,a}/\cM_{X,a}^{n+1}]^*)$$
where $[\cdots]^*$ denotes complex vector space duals.

Instead of a meromorphic map $\tau_n:Z\to\cD(\Delta_X^{(n)}/X)$, we are now interested in producing a natural meromorphic map $\mu_n$;
$$\xymatrix{
Z\ar[rr]^{\mu_n\quad\quad}\ar[dr]_\pi && \grass_{r_n}(\Diff_X^{(n)})\ar[dl]\\
& X}$$
For $a\in X$ smooth, the idea is to view $\pi^{-1}(a)$ as parametrising a family of linear subspaces of $[\cO_{X,a}/\cM_{X,a}^{n+1}]^*$.
Given $c\in\pi^{-1}(a)$, $Z_{\rho(c)}$ is an analytic set in $X$ passing through $a$.
Hence $[\cO_{Z_{\rho(c)},a}/\cM_{Z_{\rho(c)},a}^{n+1}]^*$ is a subspace of $[\cO_{X,a}/\cM_{X,a}^{n+1}]^*$.
For $a$ and $c$ sufficiently general the codimension of this subspace is constant, say $r_n$.
So $[\cO_{Z_{\rho(c)},a}/\cM_{Z_{\rho(c)},a}^{n+1}]^*$ corresponds to a point in $\grass_{r_n}([\cO_{X,a}/\cM_{X,a}^{n+1}]^*)=\grass_{r_n}(\Diff_X^{(n)})_a$.
The map $\mu_n$ will associate to $c$ this point.

One must ensure that the partial map $\mu_n$ obtained in this way is indeed meromorphic.
As we are unable to add anything enlightening on this point to what can be found in~\cite{campana}, we leave the reader to look there for the formal construction of $\mu_n$.

The (dual to the) argument for Claim~\ref{claim1} shows that if $A,B\subset X$ are irreducible analytic sets with $A\cap B\neq\emptyset$, then $A=B$ if and only if for some $a\in A\cap B$ and all $n>0$, $[\cO_{A,a}/\cM_{A,a}^{n+1}]^*=[\cO_{B,a}/\cM_{B,a}^{n+1}]^*$ as subspaces of $[\cO_{X,a}/\cM_{X,a}^{n+1}]^*$.
Then, using $\omega$-compactness as before, one shows that for sufficiently large $n$, $\mu_n$ is a bimeromorphism with its image.
As $\grass_{r_n}(\Diff_X^{(n)})\to X$ is projective (by the Pl\"ucker embedding, see~\cite{groth2}), this shows that $\pi:Z\to X$ is Moishezon, hence giving an alternative proof for Theorem~\ref{cf}.
\qed


\bigskip
\begin{thebibliography}{1}

\bibitem{campana}
F.~Campana.
\newblock Alg\'ebricit\'e et compacit\'e dans l'espace des cycles d'un espace
  analytique complexe.
\newblock {\em Mathematische Annalen}, 251(1):7--18, 1980.

\bibitem{fujiki}
A.~Fujiki.
\newblock On the {D}ouady space of a compact complex space in the category
  {$\mathcal{C}$}.
\newblock {\em Nagoya Mathematical Journal}, 85:189--211, 1982.

\bibitem{groth2}
A.~Grothendieck.
\newblock Techniques de construction en g\'eom\'etrie analytique {V}.
  {F}ibr\'es vectoriels, fibr\'es projectifs, fibr\'es en drapeaux.
\newblock {\em S\'eminaire Henri Cartan}, 13(12), 1960/61.

\bibitem{groth}
A.~Grothendieck.
\newblock Techniques de construction en g\'eom\'etrie analytique {VII}.
  \'{E}tude locale des morphisme: \'el\'ements de calcul infinit\'esimal.
\newblock {\em S\'eminaire Henri Cartan}, 13(14), 1960/61.

\bibitem{kantor}
J.M. Kantor.
\newblock Formes et op\'erateurs diff\'erentiels sur les espaces analytiques
  complexes.
\newblock {\em Bull. Soc. Math. France M\'em.}, 53:5--80, 1977.

\bibitem{ret}
R.~Moosa.
\newblock A nonstandard {R}iemann existence theorem.
\newblock {\em Transactions of the American Mathematical Society},
  356(5):1781--1797, 2004.

\bibitem{pillay}
A.~Pillay.
\newblock Model-theoretic consequences of a theorem of {C}ampana and {F}ujiki.
\newblock {\em Fundamenta Mathematicae}, 174(2):187--192, 2002.

\end{thebibliography}
\end{document}